\documentclass[11pt]{amsart}
\usepackage{amsmath,amsfonts}

\newtheorem{thm}{Theorem}[section]
\newtheorem{cor}[thm]{Corollary}
\newtheorem{lem}[thm]{Lemma}

\newtheorem{prop}[thm]{Proposition}

\newtheorem{notation}[thm]{Notation}
\numberwithin{equation}{section}

\begin{document}

\title[On diffeomorphisms deleting weakly compacta]{On diffeomorphisms deleting
weakly compacta in Banach spaces}
\author{Daniel Azagra and Alejandro Montesinos}
\date{August 21, 2002}
\thanks{D. Azagra was supported by a Marie Curie
Fellowship of the European Community Training and Mobility of
Researchers Programme under contract number HPMF-CT-2001-01175.
A. Montesinos was supported by a FPI Fellowship from the Comunidad
de Madrid}

\begin{abstract}
We prove that if $X$ is an infinite-dimensional Banach space with
$C^p$ smooth partitions of the unity then $X$ and $X\setminus K$
are $C^p$ diffeomorphic, for every weakly compact set $K\subset
X$.
\end{abstract}

\maketitle

\section[Introduction]{Introduction, main results and preliminaries}

A subset $K$ of $X$ is said to be topologically negligible
provided there exists a homeomorphism $h:X\to X\setminus K$. The
homeomorphism $h$ is usually required to be the identity outside
a given neighborhood of $K$. Here $X$ can be a Banach space, a
manifold, or just a topological space, though in this paper we
will only consider the case when $X$ is an infinite-dimensional
Banach space and $h$ is a diffeomorphism (of course, points are
{\em not} topologically negligible in finite-dimensional spaces).

The theory of topological negligibility started in 1953 when
Victor L. Klee \cite{K} proved that, if $X$ is a non-reflexive
Banach space or an infinite-dimensional $L^{p}$ space and $K$ is a
compact subset of $X$, there exists a homeomorphism between $X$
and $X\setminus K$ which is the identity outside a given
neighborhood of $K$. Klee also proved that for those
infinite-dimensional Banach spaces X the unit sphere and the unit
ball are homeomorphic to any of the closed hyperplanes in X, and
gave a topological classification of convex bodies in Hilbert
spaces. In subsequent papers, Bessaga and Klee generalized those
results to every infinite-dimensional normed space \cite{BK1,
BK2, BP}.

Klee's original proofs were of a strong geometrical flavor: very
beautiful, but rather difficult to handle in an analytical way.
Nevertheless, C. Bessaga found elegant explicit formulas for
deleting homeomorphisms, based on the existence of continuous
noncomplete (nonequivalent) norms in every infinite-dimensional
Banach space. This discovery allowed him in 1966 to construct
diffeomorphisms which delete points in the Hilbert space, and to
prove that the Hilbert space is diffeomorphic to its unit sphere
\cite{Be1}. These striking results of Bessaga's have been highly
celebrated and they remain a key ingredient in the proofs of the
already classic fundamental theorems on Hilbert manifolds (e.g.,
that every two homotopic Hilbert manifolds are diffeomorphic, see
\cite{BurKui, EE, Moulis}). Later on (and by using Bessaga's
result), J. West produced the following theorem (which remains
the most powerful result about topological negligibility in
Hilbert manifolds): if $K$ is a closed, locally compact subset of
a Hilbert manifold $M$, $U$ is an open neighborhood of $K$ and
$G$ is an open covering of $M$, then there exists a $C^\infty$
diffeomorphism $h:M\to M\setminus K$ which is the identity off
$U$ and is limited by $G$ (this roughly means that $h$ is
arbitrarily close to the identity mapping).

Apart from the classification of Hilbert manifolds by homotopy
type, the results about topological negligibility have found many
interesting applications in several branches of mathematics,
which include fixed point theory, smooth topological
classification of convex bodies, strange phenomena concerning
ordinary differential equations and dynamical systems in infinite
dimensions, the failure of Rolle's theorem in infinite dimensions
and many more things, see \cite{Ado1, Ado2, Be2, Ade2, AGJ, AJ1,
AJ2, AJ3} and the references therein. Very recently, Manuel
Cepedello and the first-named author have used smooth topological
negligibility to prove the following approximate strong version
of the Morse-Sard theorem: the smooth functions {\em with no
critical points} are dense in the space of continuous functions
on every Hilbert manifold. More precisely, if $M$ is a smooth
manifold modeled on the separable Hilbert space and
$f:M\to\mathbb{R}^{m}$ and $\varepsilon:M\to (0,\infty)$ are
continuous functions then there exists a $C^\infty$ smooth
function $g:M\to\mathbb{R}^{m}$ with no critical points and such
that $\|f-g\|\leq\varepsilon$ (a positive consequence of this
theorem is the following fact, which may be regarded as a
nonlinear analogue of the Hahn-Banach theorem: if $C_1$ and $C_2$
are two disjoint closed sets in a Hilbert manifold $M$, then
there exists a smooth submanifold of codimension one which
separates $C_1$ and $C_2$ and which is a level set of a smooth
function with no critical points on $M$). See \cite{AC2}.

In view of the interest of such applications, it is natural to
try to extend these results to Banach spaces other than the
Hilbert space.

The real-analytic and smooth negligibility of compact sets in
Banach spaces was studied by Tadeusz Dobrowolski \cite{Do}, who
developed Bessaga's non-complete norm technique in the smooth
case and generalized some of the results of \cite{Be1, BP}. He
\cite{Do} showed that for every infinite-dimensional Banach space
$X$ having a $C^{p}$ non-complete norm, and for every compact set
$K$ in $X$, the space $X$ is $C^{p}$ diffeomorphic to $X\setminus
K$. Unfortunately, it is still unknown whether every Banach space
with a $C^p$ smooth equivalent norm possesses a noncomplete $C^p$
smooth norm as well.

Without showing the existence of smooth non-complete norms, the
first-named author proved \cite{A, thesis} that every Banach space
$(X,\|\cdot\|)$ with a $C^{p}$ smooth norm $\varrho$ is $C^{p}$
diffeomorphic to $X\setminus\{0\}$ and, moreover, that every
hyperplane $H$ in $X$ is $C^{p}$ diffeomorphic to the sphere
$\{x\in X\mid \varrho(x)=1\}$.

In subsequent work \cite{Ado1, Ado2}, T. Dobrowolski and the
first-named author strengthened the new technique of deleting
points introduced in \cite{A} so as to generalize some results and
applications of smooth negligibility of compacta and subspaces to
the class of all Banach spaces having $C^{p}$ smooth norms.

Despite all these efforts, the natural question as to the
characterization of those Banach spaces in which compact sets are
topologically negligible remains open. This is due to a surprising
(and rather uncomfortable) theorem proved by R. Haydon
\cite{Haydon1, Haydon2}: there are Banach spaces which have
$C^\infty$ smooth bump functions, and even $C^\infty$ smooth
partitions of unity, but do not possess any equivalent $C^1$
smooth norm.

In this paper we deal with the following natural question: what
can be said about smooth negligibility of compacta in those
Banach spaces with smooth partitions of unity? As we have just
pointed out, there are Banach spaces with smooth partitions of
unity which have no equivalent smooth norms, and therefore the
known results on diffeomorphisms deleting compacta are useless in
this setting. Nevertheless, we will prove in this paper that every
(weakly) compact subset $K$ of an infinite-dimensional Banach
space with $C^p$ smooth partitions of unity can be removed by a
$C^p$ diffeomorphism $h:X\to X\setminus K$ which is the identity
outside a given starlike body $A$ containing $K$. More precisely,
our main result is as follows.

\begin{thm}\label{main thm}
Let $X$ be an infinite-dimensional Banach space which has $C^p$
smooth partitions of the unity. Then, for every weakly compact set
$K\subset X$ and every starlike body $A$ such that
$\textrm{dist}(K,X\setminus A)>0$, there exists a $C^p$
diffeomorphism $h:X\longrightarrow X\setminus K$ such that $h$ is
the identity outside $A$.

In particular, when $K$ is compact and $K\subset\textrm{int}(A)$,
there always exists such a deleting diffeomorphism $h$.
\end{thm}
Here, $p\in\mathbb{N}\cup\{\infty\}$. Of course if $X$ is
finite-dimensional, there exist no such deleting diffeomorphisms.

The class of Banach spaces which admit smooth partitions of unity
is quite large. For instance, every Banach space with a separable
dual admits $C^1$ smooth partitions of unity, so does every
reflexive space, and there are many other simple conditions that
ensure the existence of smooth partitions of unity in a Banach
space; see \cite{DGZ}. On the other hand, it is an open problem
to know whether every Banach space with a $C^p$ smooth equivalent
norm has $C^p$ smooth partitions of unity (see \cite{DGZ}). If a
positive answer to this question is ever reached, then Theorem
\ref{main thm} will be an extension of the main theorem in
\cite{Ado1}. Otherwise and for the time being, by combining
Theorem \ref{main thm} with the main result of \cite{Ado1}, we
may deduce the following.

\begin{cor}
Let $X$ be an infinite-dimensional Banach space. Assume that
either $X$ possesses a $C^p$ smooth norm or else $X$ has $C^p$
smooth partitions of unity. Then, for every compact set $K\subset
X$ and every $C^p$ smooth starlike body $A$ such that
$K\subset\textrm{int}(A)$, there exists a $C^p$ diffeomorphism
$h:X\longrightarrow X\setminus K$ such that $h$ is the identity
outside $A$.
\end{cor}
It should be noted that, for the time being, no one knows of an
infinite-dimensional Banach space with a $C^1$ bump function
which does not have either a $C^1$ smooth norm or $C^1$ smooth
partitions of unity (hence which does not fall into the category
to which the above Corollary applies). On the other hand, it is
easy to see that the existence of a $C^1$ smooth bump is a
necessary condition for a Banach space $X$ to have a
diffeomorphism from $X$ onto $X\setminus\{0\}$ which restricts to
the identity outside some ball.

We should also stress that, as a consequence of Theorem \ref{main
thm}, all of the corollaries of the main theorem of \cite{Ado1}
proved in that paper are true for every infinite-dimensional
Banach space having $C^p$ smooth partitions of unity. For
instance, Garay's results concerning strange phenomena for ODEs in
Banach spaces can be readily extended to this category. We will
not elaborate on these topics; see \cite{Ado1} and the references
therein.

At this point we need to introduce some terminology and notation
concerning starlike bodies, which, apart from the statements of
the preceding results, will play a key role in our proofs. The
number of things that can be proved if one only knows that there
are {\em enough} smooth starlike bodies in our space is somewhat
surprising. In fact, the requirement that our space $X$ has $C^p$
smooth partitions of unity will only be used to ensure that our
space has enough smooth starlike bodies.

A closed subset $A$ of a Banach space $X$ is said to be a
starlike body if there exists a point $a_{0}$ in the interior of
$A$ such that every ray emanating from $a_{0}$ meets $\partial
A$, the boundary of $A$, at most once. We will say that $a_{0}$
is a center of $A$. There can obviously exist many centers for a
given starlike body.  Up to a suitable translation, we can always
assume that $a_{0}=0$ is the origin of $X$, and we will often do
so, unless otherwise stated. For a starlike body $A$ with center
$a_{0}$, we define the characteristic cone of $A$ as $$ cc
A=\{x\in X| a_{0}+r(x-a_{0})\in A\ \text{for\  all}\ r>0\}, $$
and the Minkowski functional of $A$ with respect to the center
$a_{0}$ as $$ \mu_{A, a_{0}}(x)=\mu_{A}(x)=\inf\{t>0 \, | \,\,
x-a_{0}\in t(-a_{0}+A)\}\, \textrm{ for all $x\in X$.} $$ Note
that $\mu_{A}(x)=\mu_{-a_{0}+A}(x-a_{0}) \textrm{ for all } x\in
X.$ It is easily seen that $\mu_{A}$ is a continuous function
which satisfies $\mu_{A}(a_{0}+rx)=r\mu_{A}(a_{0}+x)$ for every
$r\geq 0$ and $x\in X$, and $\mu_{A}^{-1}(0)=cc A$. Moreover,
$A=\{x\in X| \mu_{A}(x)\leq 1\}$, and $\partial A=\{x\in X \mid
\mu_{A}(x)=1\}$. Conversely, if $\psi:X\to [0, \infty)$ is
continuous and satisfies $\psi(a_{0}+\lambda
x)=\lambda\psi(a_{0}+x)$ for all $\lambda\geq 0$, then
$A_{\psi}=\{x\in X \mid \psi(x)\leq 1\}$ is a starlike body. More
generally, for a continuous function $\psi:X\to[0,\infty)$ such
that $\psi_x(\lambda)=\psi(a_{0}+\lambda x)$, $\lambda>0$, is
increasing and $\sup\{\psi_x(\lambda)
: \lambda>0\}>\varepsilon$ for every $x\in X\setminus \psi^{-1}(0)$,
the set $\psi^{-1}([0,\varepsilon])$ is a starlike body whose
characteristic cone is $\psi^{-1}(0)\ni a_{0}$.

A familiar important class of starlike bodies are {\em convex
bodies}, that is, starlike bodies that are convex. For a convex
body $U$, $cc U$ is always a convex set, but in general the
characteristic cone of a starlike body is not convex. Starlike
bodies can also be related to $n$-homogeneous polynomials, since
the level sets of such polynomials are always boundaries of
starlike bodies.

We will say that $A$ is a $C^{p}$ smooth starlike body provided
its Minkowski functional $\mu_{A}$ is $C^{p}$ smooth on the set
$X\setminus cc A =X\setminus \mu_{A}^{-1}(0)$. This is equivalent
to saying that $\partial A$ is a $C^p$ smooth one-codimensional
submanifold of $X$ such that no affine hyperplane tangent to
$\partial A$ contains a ray emanating from the center $a_{0}$.
Throughout this paper, $p=0,1,2, ...., \infty$, and $C^0$ smooth
means just continuous.

We will also say that $A$ is Lipschitz if $\mu_{A}$ is a
Lipschitz function on $X$. It is easy to see that every {\em
convex} body is Lipschitz with respect to any point in its
interior (but this is no longer true if we drop convexity: even
in the plane $\mathbb{R}^{2}$ there are starlike bodies which are
not Lipschitz).

All the starlike bodies that we will deal with in this paper are
{\em radially bounded}. A starlike body $A$ is said to be radially
bounded provided that, for every ray emanating from the center
$a_{0}$ of $A$, the intersection of this ray with $A$ is a
bounded set. This amounts to saying that $cc A=\{a_{0}\}$.

In finite dimensions every radially bounded starlike body is in
fact bounded (because the Minkowski functional of the body
attains an absolute minimum on the unit sphere, which is compact),
but this is no longer true in infinite-dimensional Banach spaces.
For instance, $A=\{x\in\ell_{2} :
\sum_{n=1}^{\infty}x_{n}^{2}/2^{n}\leq 1\}$ is a radially bounded
convex body which is not bounded in the Hilbert space $\ell_{2}$;
the body $A$ is the unit ball of the nonequivalent $C^\infty$
smooth norm $\omega(x)=\sum_{n=1}^{\infty}x_{n}^{2}/2^{n}$ in
$\ell_{2}$.

For every {\em bounded} starlike body $A$ in a Banach space $(X,
\|\cdot\|)$ there are constants $M, m>0$ such that
    $$
    m\|x\|\leq \mu_{A}(x)\leq M\|x\| \textrm{ for all } x\in X.
    $$
If $A$ is just {\em radially bounded} then we can only ensure that
    $$
    \mu_{A}(x)\leq M\|x\| \textrm{ for all } x\in X,
    $$
for some $M>0$. As is shown implicitly in \cite[Proposition
II.5.1]{DGZ}, a Banach space $X$ has a $C^p$ smooth bump function
if and only if there is a bounded $C^p$ smooth starlike body in
$X$. The reader might want to consult the references \cite{AC1,
Ade2, Ado3, AJ2, AJ3} for other properties of starlike bodies.

We will finish these preliminaries with some nonstandard notation
concerning strict inclusions between starlike bodies. In our
proofs we will often require that, for a couple of starlike bodies
$A\subset B$, the boundaries of $A$ and $B$ are well separated.
There are at least two nonequivalent natural notions of
separation between boundaries of starlike bodies, and we will
need to use both of them, as each one has its own advantages. The
strongest and most natural notion corresponds to the fact that
the distance between $A$ and $X\setminus B$ is positive. We will
use the notation
    $$
    A\subset_{d} B
    $$
to mean that $\textrm{dist}(A,X\setminus B)>0$, and we will say
that $B$ strictly contains $A$ in the distance sense. Notice that
this notion makes sense even though $A$ and $B$ do not have the
same center, or even if $A$ and $B$ are mere sets, not
necessarily starlike.

The other useful notion is that the Minkowski functionals of $A$
and $B$ are well separated, in the following sense. First, note
that if $A\subseteq B$ are starlike with respect to the same
center $a_{0}$ then we always have that
$\mu_{B}(x)\leq\mu_{A}(x)$ for all $x\in X$. If we also know that
    $$
    \sup_{x\in A}\mu_{B}(x) <1
    $$
then we will denote $A\subset_{\mu}B$, saying that $B$ strictly
contains $A$ in the gauge sense. This is equivalent to saying
that there exists some $\delta>0$ such that
    $$
    a_{0}+(1+\delta)(-a_{0}+A)\subseteq B.
    $$
Of course, this notion only makes sense when $A$ and $B$ have at
least one center $a_{0}$ in common. It is immediate to see that
$A\subset_{d} B$ implies that $A\subset_{\mu} B$. The converse is
false in general, unless $A$ is Lipschitz. When $A\subset B$ have
the same center and $A$ is Lipschitz we have that $A\subset_{d} B$
if and only if $A\subset_{\mu} B$ (see Lemma \ref{Lipschitz
implies that mu implies d} below).


\medskip

\section[Proof]{Proof of the main result}

\medskip

In contrast with Bessaga-type constructions \cite{Be1, Do, A,
thesis, Ado1, Ado2}, our proof does not provide an explicit
elegant formula for the deleting diffeomorphism. The reason why
we cannot use those Bessaga-type formulas in this setting is the
following one. Such deleting formulas are of the form
$h:X\setminus K\to X, h(x)=x+p(f(x))$, where $p:(0,\infty)\to X$
is a {\em deleting path}, $f$ is a function such that
$f^{-1}(0)=K$, and both $f$ and $p$ satisfy a Lipschitz condition
with respect to the Minkowski functional $\omega$ of a convex body
which is radially bounded but not bounded (for instance $\omega$
could be a continuous noncomplete norm). The path $p$ can always
be assumed to be $C^\infty$ smooth, but the function $f$ cannot,
in general. One could think that if one approximates the function
$f$ well enough by a smooth function $g$ then the formula $x\to
x+p(g(x))$ should define a diffeomorphism from $X\setminus K$ onto
$K$. This approach can only be successful if we further ensure
that $g$ is still Lipschitz with respect to $\omega$.
Unfortunately, for an infinite-dimensional Banach space $X$ with
smooth partitions of unity, it is unknown whether a given
function $f$ which is Lipschitz with respect to a continuous norm
$\omega$ can be uniformly approximated by smooth functions which
are still Lipschitz with respect to $\omega$; in fact the
question is open even when the norm $\omega$ is complete.

So we will rather turn to the origins and find inspiration in the
geometrical ideas of the pioneering work of Klee's \cite{K}. We
will need to consider an infinite composition of carefully
constructed self-diffeomorphisms of $X$.

The main ingredient of our proof is the following Proposition,
which implies that if our infinite-dimensional space $X$ has
enough smooth starlike bodies then every weakly compact set $K$
can be removed by means of a diffeomorphism $h:X\to X\setminus K$
which is the identity outside some starlike body.

\begin{prop}\label{construction of the deleting diffeomorphism}
Let $X$ be a Banach space, and $K$ a subset of $X$. Assume that
there are sequences $(P_{n})$, $(C_{n})$, $(A_{n})$, $(B_{n})$,
$(Q_{n})$, $(D_{n})$, $(E_{n})$ of subsets of $X$ and a sequence
$(c_{n})$ of points of $X$ satisfying the following conditions for
each $n\in\mathbb{N}$:
\begin{enumerate}
\item $A_{n}$, $B_{n}$, $Q_{n}$, $D_{n}$, $E_{n}$ are radially
bounded $C^p$ smooth starlike bodies with respect to $c_{n+2}$;
\item $C_{n+2}\subset D_{n}\subset_{\mu} E_{n}\subset_{\mu} A_{n}\subset C_{n+1}
\subset P_{n+1}\subset B_{n}\subset_{\mu} Q_{n}\subset P_{n}$
\item $\bigcap_{n=1}^{\infty}C_{n}=\emptyset$
\item $\bigcap_{n=1}^{\infty}P_{n}=K$.
\end{enumerate}
Then there exists a $C^p$ diffeomorphism $\Psi:X\longrightarrow
X\setminus K$ such that $\Psi$ is the identity on $X\setminus
P_{1}$.
\end{prop}
In order to prove this Proposition we will only require a simple
geometrical Lemma.
\begin{lem}[The four bodies lemma]\label{four bodies}
Let $X$ be a Banach space, and let $A, B, C, D$ be four radially
bounded $C^p$ smooth starlike bodies with respect to the same
point $a_{0}\in\textrm{int}(A)$. Assume that
    $$
    A\subset_{\mu} B\subset C\subset_{\mu} D.
    $$
Then there exists a $C^p$ diffeomorphism $h:X\to X$ such that
\begin{enumerate}
\item $h(B)=C$
\item $h$ is the identity on $A\cup(X\setminus D)$.
\end{enumerate}
\end{lem}
\begin{proof}
We may assume $a_{0}=0$. Since $A\subset_{\mu} B$ and
$C\subset_{\mu} D$, there exists some $\delta\in (0,1)$ such that
$A\subset (1-\delta) B$ and $(1+\delta)C\subset D$. Take a
$C^\infty$ smooth function $\lambda:\mathbb{R}\to\mathbb{R}$ such
that $\lambda$ is non-decreasing, $\lambda(t)=0$ if $t\leq
1-\delta$, and $\lambda(t)=1$ for $t\geq 1$. Define then $f:X\to
X$ by
    $$
    f(x)=\big[\lambda(\mu_{B}(x))\frac{\mu_{B}(x)}{\mu_{C}(x)}+
    1-\lambda(\mu_{B}(x))\big] x, \, \textrm{ if } x\neq 0,
    $$
and $f(0)=0$. It is easy to check that $f$ is a $C^p$
diffeomorphism of $X$ such that $f(B)=C$ and $f$ is the identity
on $A$.

On the other hand, pick $\theta:\mathbb{R}\to\mathbb{R}$ a
$C^\infty$ smooth function such that $\theta$ is non-increasing,
$\theta(t)=1$ if $t\leq 1+\delta/4$, and $\theta(t)=0$ if $t\geq
1+\delta/2$. Consider the mapping $g:X\setminus\{0\}\to
X\setminus\{0\}$ defined by
    $$
    g(x)=\big[\theta(\mu_{C}(x))\frac{\mu_{C}(x)}{\mu_{B}(x)}+
    1-\theta(\mu_{C}(x))\big] x,
    $$
which is a $C^p$ diffeomorphism as well. Now define $h:X\to X$ by
    $$
    h(x)=\left\{
    \begin{array}{ll}
    &f(x)\quad \hspace{4mm} \textrm{if $\mu_{B}(x)<1+\frac{\delta}{4}$;}\\
    &g^{-1}(x) \quad \textrm{if $1<\mu_{B}(x)$}
    \end{array}\right.
    $$
Observe that if $1\leq\mu_{B}(x)\leq 1+\delta/4$ then
$f(x)=\big[\mu_{B}(x)/\mu_{C}(x)\big] x=g^{-1}(x)$ ; hence $h$ is
well-defined and locally a $C^p$ diffeomorphism. Moreover, it is
easy to see that $h(X\setminus(1+\delta/4)B)=X\setminus
(1+\delta/4)C$, which (bearing in mind the definition of $h$)
implies that $h$ is one-to-one. On the other hand, since
$h((1+\delta/4)B)=(1+\delta/4)C$ and $h(X\setminus
B)=g^{-1}(X\setminus B)=X\setminus C$, it follows that $h$ is a
surjection. Therefore $h:X\to X$ is a $C^p$ diffeomorphism.
Finally, it is clear that $h(B)=C$, and $h$ is the identity on
$A\cup \big(X\setminus (1+\delta/2)B\big)\supset
A\cup\big(X\setminus D\big)$.
\end{proof}

\medskip
\begin{center}
{\bf Proof of Proposition \ref{construction of the deleting
diffeomorphism}}
\end{center}

\noindent Fix any $n\in\mathbb{N}$. Consider the inclusions of
bodies
\begin{eqnarray*}
& &D_{n}\subset_{\mu} E_{n}\subset_{\mu}B_{n}\subset_{\mu}Q_{n}\\
& &D_{n}\subset_{\mu} E_{n}\subset_{\mu}A_{n}\subset_{\mu}Q_{n}.
\end{eqnarray*}
According to the Four Bodies Lemma there exist $C^p$
diffeomorphisms $f_{n}, g_{n}:X\to X$ such that
\begin{eqnarray*}
& &f_{n}(E_{n})=B_{n}, \textrm{ and $f_{n}$ is the identity on
$D_{n}\cup(X\setminus Q_{n})$},\\ & &g_{n}(E_{n})=A_{n},
\textrm{and $g_{n}$ is the identity on $D_{n}\cup(X\setminus
Q_{n})$}.
\end{eqnarray*}
Define then $h_{n}=g_{n}\circ f_{n}^{-1}:X\to X$, which is a $C^p$
diffeomorphism of $X$ satisfying that
    $$
    h_{n}(B_{n})=A_{n}, \textrm{ and }
    h_{n} \textrm{ is the identity on } D_{n}\cup(X\setminus
    Q_{n}).
    $$
Now consider the family of $C^p$ diffeomorphisms $(h_{n})$. For
each $n\in\mathbb{N}$ define the mapping $\psi_{n}:X\to X$ by the
composition
    $$
    \psi_{n}(x)=(h_{1}\circ h_{2}\circ ...
    \circ h_{n-1}\circ h_{n})(x),
    $$
which is obviously a $C^p$ diffeomorphism of $X$. Since $h_{n}$
is the identity on $X\setminus Q_{n}$ and $Q_{n}\subset P_{n}$,
we have that $h_{n}$ is the identity on $X\setminus P_{n}$. It
follows that
    $$
    {\psi_{n}}_{|X\setminus P_{n}}={\psi_{n-1}}_{|X\setminus P_{n}}
    \textrm{ for all } n\geq 2. \eqno(1)
    $$
Note that, from the conditions in the statement of Proposition
\ref{construction of the deleting diffeomorphism}, we know that
$$X\setminus P_{n}\subset X\setminus P_{n+1}\subset X\setminus K,
\textrm{ for all $n$, and } X\setminus
K=\bigcup_{n=1}^{\infty}X\setminus P_{n}.\eqno(2)$$ Then we can
define $\psi:X\setminus K\to X$ by letting
    $$
    \psi_{|X\setminus P_{n+1}}={\psi_{n}}_{|X\setminus P_{n+1}}
    \textrm{ for each } n\in\mathbb{N}. \eqno(3)
    $$
Taking equations $(1)$ and $(2)$ above into account, it is clear
that the mapping $\psi$ is well defined, one-to-one, and is
locally a $C^p$ diffeomorphism. Let us see that $\psi$ is
surjective and therefore a $C^p$ diffeomorphism from $X\setminus
K$ onto $X$.

Bearing in mind that $h_{j}$ is the identity on $D_{j}\supset
C_{j+2}$ and $A_{j}\subset C_{j+1}$, we have that
$h_{j}(A_{n})=A_{n}$ if $j\leq n-1$, and since
$h_{n}(B_{n})=A_{n}$ we may deduce that
    $$
    \psi_{n}(B_{n})=h_{1}\circ ...\circ h_{n}(B_{n})=
    h_{1}\circ...\circ h_{n-1}(A_{n})=h_{1}\circ...\circ h_{n-2}(A_{n})
    =...=A_{n};
    $$
and in particular $\psi_{n}(X\setminus B_{n})=X\setminus A_{n}$.
But, by the hypothesis on the bodies, $P_{n+1}\subset B_{n}\subset
P_{n}$, that is $X\setminus P_{n}\subset X\setminus B_{n}\subset
X\setminus P_{n+1}$, and hence $$\psi(X\setminus
B_{n})=\psi_{n}(X\setminus B_{n})=X\setminus A_{n}. \eqno (4)$$
Now, note that the hypothesis of Proposition \ref{construction of
the deleting diffeomorphism} imply that $C_{n+2}\subset
A_{n}\subset C_{n+1}$, $\bigcap_{n=1}^{\infty}C_{n}=\emptyset$,
which yield
    $$
    X=\bigcup_{n=1}^{\infty}(X\setminus A_{n}). \eqno(5)
    $$
On the other hand, since $K=\bigcap_{n=1}^{\infty}P_{n+1}\subset
\bigcap_{n=1}^{\infty}B_{n}\subset
\bigcap_{n=1}^{\infty}P_{n}=K$, we have that
    $$
    X\setminus K=\bigcup_{n=1}^{\infty}(X\setminus B_{n}).
    \eqno(6)
    $$
Now, by combining equations $(4)$, $(5)$ and $(6)$, we get that
    $$
    \psi(X\setminus K)=\psi\big(\bigcup_{n=1}^{\infty}(X\setminus B_{n})\big)=
    \bigcup_{n=1}^{\infty}(X\setminus A_{n})=X,
    $$
hence $\psi$ is a $C^p$ diffeomorphism from $X\setminus K$ onto
$X$. Moreover, if $x\in X\setminus P_{1}\subset X\setminus P_{2}$,
from the definition of $\psi$, and bearing in mind that $h_{1}$
is the identity on $X\setminus P_{1}$, we conclude that
$\psi(x)=\psi_{1}(x)=h_{1}(x)=x$. Finally, if we define
$\Psi=\psi^{-1}$, it is clear that $\Psi$ is a $C^p$
diffeomorphism from $X$ onto $X\setminus K$ which is the identity
off $P_{1}$.

\medskip

The next step in the proof of our main theorem is of course to
ensure that if an infinite-dimensional Banach space $X$ has $C^p$
smooth partitions of unity then, for every weakly compact set
$K\subset X$, there are families of $C^p$ smooth starlike bodies
satisfying the conditions of Proposition \ref{construction of the
deleting diffeomorphism}.

\begin{prop}\label{existence of enough starlike bodies}
Let $X$ be an infinite-dimensional Banach space which admits $C^p$
smooth partitions of unity. Then there exists $B$, a radially
bounded $C^p$ smooth starlike body with respect to the origin,
such that, for every weakly compact set $K\subset X$ and for
every bounded starlike body $A\supset K$, there is some $r>0$ such
that $A\subset rB$, and there are sequences $(P_{n})$, $(C_{n})$,
$(A_{n})$, $(B_{n})$, $(Q_{n})$, $(D_{n})$, $(E_{n})$ of subsets
of $X$ and a sequence $(c_{n})$ of points of $X$ satisfying the
following conditions for each $n\in\mathbb{N}$:
\begin{enumerate}
\item $A_{n}$, $B_{n}$, $Q_{n}$, $D_{n}$, $E_{n}$ are radially
bounded $C^p$ smooth starlike bodies with respect to $c_{n+2}$;
\item $C_{n+2}\subset D_{n}\subset_{\mu} E_{n}\subset_{\mu} A_{n}\subset C_{n+1}
\subset P_{n+1}\subset B_{n}\subset_{\mu} Q_{n}\subset P_{n}$;
\item $\bigcap_{n=1}^{\infty}C_{n}=\emptyset$;
\item $\bigcap_{n=1}^{\infty}P_{n}=K$;
\item $P_{1}\subset rB$.
\end{enumerate}
\end{prop}
In the sequel such a body $B$ will be called {\em universal}.

The proof of Proposition \ref{existence of enough starlike
bodies} is quite long and will be split into several lemmas.
\begin{notation}
{\em If $X$ is a Banach space and $B_{X}=\{x\in X : \|x\|\leq 1\}$
is its unit ball, for all subsets $A, B$ of $X$ and for every
$\varepsilon >0$, we will denote
    $$
    [A,B]=\{tx+(1-t)y : x\in A, y\in B, t\in [0,1]\},
    $$
and $N(A,\varepsilon)=\{x\in X : \textrm{dist}(x,
A)\leq\varepsilon\}=\overline{A+\varepsilon B_{X}}$. When
$A=\{a\}$ is a singleton we will simply write $[A,B]=[a,B]$.}
\end{notation}

\begin{lem}\label{that union of cones is starlike}
Let $X$ be a Banach space, $C$ a bounded convex body in $X$, and
$K$ a weakly compact subset of $X$. Then $V:=[K, C]$ is a starlike
body with respect to every interior point of $C$. Moreover, $V$ is
bounded and $\mu_{V}:X\to[0,\infty)$ is Lipschitz.
\end{lem}
\begin{proof}
Since $C\subseteq V$, it is obvious that $V$ has nonempty
interior.

Let us see that $V$ is closed. Take a sequence $(z_{n})$ in $V$
converging to a point $z_{0}$. Each $z_{n}$ is of the form
$z_{n}=t_{n}x_{n}+(1-t_{n})y_{n}$, with $t_{n}\in [0,1]$,
$x_{n}\in C$, $y_{n}\in K$. Since $K$ is weakly compact and
$[0,1]$ is compact, we may assume, passing to a subsequence if
necessary, that $t_{n}\to t_{0}\in [0,1]$ and $y_{n}\to y_{0}\in
K$ weakly. Then we distinguish two possibilities: either
$t_{0}\neq 0$ or $t_{0}=0$. If $t_{0}\neq 0$ then we see that
$x_{n}=t_{n}^{-1}(z_{n}-(1-t_{n})y_{n})$ weakly converges to the
point $x_{0}:=t_{0}^{-1}(z_{0}-(1-t_{0})y_{0})$, which must
belong to $C$ because $C$ is closed and convex, hence weakly
closed. Then $z_{0}=\lim_{n}z_{n}=w\textrm{-}\lim_{n}
z_{n}=t_{0}x_{0}+(1-t_{0})y_{0}$ clearly belongs to $V$. On the
other hand, if $t_{0}=0$ then, since $C$ is bounded, we have that
$\|t_{n}x_{n}\|\to 0$, and we deduce that $z_{0}=\lim_{n}
z_{n}=w\textrm{-}\lim_{n}z_{n}=w\textrm{-}\lim_{n} y_{n}=y_{0}\in
K\subset V$. In either case, $z_{0}\in V$, and this shows that
$V$ is closed.

Now let us see that $V$ is starlike with respect to every point
$x_{0}\in\textrm{int}(C)$. Take two points $x_{1},
x_{2}\in\partial V\subset V$ with $x_{1}\in [x_{0}, x_{2}]$.
Assuming that $x_{1}\neq x_{2}$ we will get a contradiction.
Indeed, since
    $$
    V=\bigcup_{y\in K}[y, C]
    $$
and $x_{1}\in\partial V$, we have that $x_{1}\in
X\setminus\textrm{int}([y,C])$ for every $y\in K$. Hence, for
every $y\in K$, either $x_{1}\in\partial[y,C]$ or $x_{1}\notin
[y,C]$; in either case, since $[y,C]$ is a starlike body with
respect to $x_{0}\in\textrm{int}(C)$, and $x_{2}\neq x_{1}\in
[x_{0}, x_{2}]$, we get that $x_{2}\notin [y,C]$. But then we have
that $$x_{2}\notin\bigcup_{y\in K}[y,C]=V,$$ a contradiction.

It is obvious that $V$ is bounded. It only remains to show that
$V$ is Lipschitz, that is, its Minkowski functional $\mu_{V}$
(with respect to any point $x_{0}\in\textrm{int}(C)$) is
Lipschitz. Without loss of generality we may assume that the
given center is $x_{0}=0$. Let $M>0$ be such that $\mu_{C}(x)\leq
M\|x\|$ for all $x\in X$. Since $C\subseteq [y,C]$ we have that
    $$
    \mu_{[y,C]}(x)\leq \mu_{C}(x)\leq M\|x\|
    $$
for all $x\in X$ and, bearing in mind that $[y,C]$ is a convex
body, this means that $\mu_{[y,C]}$ is $M$-Lipschitz for all
$y\in K$. On the other hand, it is easily seen that
    $$
    \mu_{V}(x)=\inf_{y\in K}\mu_{[y,C]}(x).
    $$
Now we can show that $\mu_{V}$ is $M$-Lipschitz as well. For any
given $\varepsilon>0$, $x,z\in X$, by using the above formula for
$\mu_{V}$ and the definition of inf, we obtain some $y\in K$ such
that
    $$
    \mu_{V}(x)-\mu_{V}(z)\leq\mu_{[y,C]}(x)+\varepsilon -\mu_{[y,C]}(z)\leq
    M\|x-z\|+\varepsilon.
    $$
This implies that $\mu_{V}(x)-\mu_{V}(z)\leq M\|x-z\|$ for all
$x,z\in X$, and therefore $\mu_{V}$ is $M$-Lipschitz.
\end{proof}

\begin{lem}\label{Lipschitz implies that mu implies d}
Let $X$ be a Banach space, $A$ a Lipschitz starlike body with
respect to the origin. Then, for every $\varepsilon>0$ there
exists $\delta>0$ so that $A+\delta B_{X}\subset
(1+\varepsilon)A$.
\end{lem}
\begin{proof}
Let $M$ be a Lipschitz constant for $\mu_{A}$. For a given
$\varepsilon>0$ choose $\delta>0$ with $\delta M<\varepsilon$.
Take $x=y+z$, with $y\in A$, $z\in\delta B_{X}$. Then we have
    $$
    \mu_{A}(x)=\mu_{A}(y+z)-\mu_{A}(y)+\mu_{A}(y)\leq
    M\|z\|+\mu_{A}(y)\leq M\delta +1<1+\varepsilon.
    $$
This shows that $A+\delta B_{X}\subset (1+\varepsilon)A$.
\end{proof}

\begin{lem}\label{smaller multiples are at positive distance}
Let $C$ a bounded convex body in a Banach space $X$, with
$0\in\textrm{int}(C)$. Then, for every $\delta\in (0,1)$,
$dist\big((1-\delta)C, X\setminus C\big)>0$, that is,
$(1-\delta)C\subset_{d} C$.
\end{lem}
\begin{proof}
This is an easy consequence of the preceding lemma and the fact
that $\mu_{C}$ is Lipschitz whenever $C$ is a convex body
(ensured in turn by Lemma \ref{that union of cones is starlike}
if we take $K=\{0\}$, $V=[0,C]=C$).
\end{proof}

\begin{lem}\label{pullback}
Let $T:X\longrightarrow Y$ be a continuous linear injection
between two Banach spaces. Then, for every radially bounded $C^p$
smooth body $B'$ in $Y$ which is starlike with respect to a point
$b'\in T(X)$, we have that $B=T^{-1}(B')$ is a radially bounded
$C^p$ smooth starlike body in $X$ with respect to $b=T^{-1}(b')$.
\end{lem}
\begin{proof}
Let $b'=T(b)$ be the center of $B'$. Then $A':=-b'+B'$ is starlike
with respect to the origin, radially bounded and $C^p$ smooth.
Consider the function $\psi:X\longrightarrow [0,\infty)$ defined
by $\psi(x)=\mu_{A'}(T(x))$. This function is continuous,
positively homogeneous, $C^p$ smooth on $X\setminus\{0\}$, and
$\psi(x)=0$ if and only if $x=0$. Therefore
    $$A:=\{x\in X : \psi(x)\leq 1\}
    $$
is a $C^p$ smooth starlike body in $X$ (with respect to the
origin); besides, since $\psi(x)>0$ whenever $x\neq 0$, we have
that $cc A=\{0\}$, that is, $A$ is radially bounded. It is
obvious that $A=T^{-1}(A')$. Then we see that
    $$
    B=T^{-1}(B')=T^{-1}(b'+A')=b+A
    $$
is a radially bounded $C^p$ smooth starlike body with respect to
$b\in X$.
\end{proof}

\begin{lem}\label{pullback preserves strict inclusions}
Let $T:X\longrightarrow Y$ be a continuous linear injection
between two Banach spaces. Assume that $A'$ and $B'$ are starlike
bodies with respect to $y_{0}=T(x_{0})\in T(X)$, and
$A'\subset_{\mu} B'$. Then
$A:=T^{-1}(A')\subset_{\mu}T^{-1}(B'):=B$.
\end{lem}
\begin{proof}
We may assume $x_{0}=0$. According to the proof of the preceding
lemma, $\mu_{A}=\mu_{A'}\circ T$, and $\mu_{B}=\mu_{B'}\circ T$.
Then, for every $x\in A$,
    $$
    \mu_{B}(x)=\mu_{B'}(T(x))\leq\sup_{y\in A'}\mu_{B'}(y)<1
    $$
because $A'\subset_{\mu} B'$, and therefore
    $$
    \sup_{x\in A}\mu_{B}(x)\leq\sup_{y\in A'}\mu_{B'}(y)<1,
    $$
which means that $A\subset_{\mu} B$.
\end{proof}
The following lemmas show how one can approximate and interpolate
starlike bodies with smooth starlike bodies, provided the space
has smooth partitions of unity.
\begin{lem}\label{smooth approximation of starlike bodies}
Let $X$ be a Banach space with $C^p$ smooth partitions of unity,
and $C$ a starlike body with $cc C=\{0\}$. Then, for every
$\delta>0$, there exists a $C^p$ smooth starlike body with $cc
A=\{0\}$ such that $(1-\delta)C\subset A\subset (1+\delta) C$.
\end{lem}
\begin{proof}
Since $X$ has $C^p$ smooth partitions of unity, it has a $C^p$
smooth bump as well, and in particular there exists $B$, a bounded
$C^p$ smooth starlike body with respect to the origin
\cite[Proposition II.5.1]{DGZ}. Choose $\varepsilon_{0}\in (0,1)$
such that
    $$
    \frac{1}{1-\varepsilon_{0}}<1+\delta, \textrm{ and }
    1+\varepsilon_{0}<\frac{1}{1-\delta}.
    $$
Define $\varepsilon:X\setminus\{0\}\to (0,\infty)$ by
    $$
    \varepsilon(x)=\varepsilon_{0}\mu_{C}(x) \textrm{ for all } x\neq 0,
    $$
which is a continuous strictly positive function. Since $X$ has
$C^p$ smooth partitions of unity, so does its open subset
$X\setminus\{0\}$, and therefore every continuous function on
$X\setminus\{0\}$ can be $\varepsilon$-approximated by a $C^p$
smooth function on $X\setminus\{0\}$. Hence, given the continuous
function $\mu_{C}:X\setminus\{0\}\to(0,\infty)$, there exists a
$C^p$ smooth function $g:X\setminus\{0\}\to\mathbb{R}$ such that
$|\mu_{C}(x)-g(x)|\leq\varepsilon(x)$ for all $x\neq 0$. Now
define $\psi:X\longrightarrow\mathbb{R}$ by
    $$
    \psi(x)=\mu_{B}(x)g\big(\frac{x}{\mu_{B}(x)}\big) \textrm{ if } x\neq 0,
    $$
and $\psi(0)=0$. The function $\psi$ is clearly continuous on $X$,
$\psi$ is of class $C^p$ on $X\setminus\{0\}$, and $\psi$ is
positively homogeneous. Moreover,
\begin{eqnarray*}
& &|\psi(x)-\mu_{C}(x)|=
\big|\mu_{B}(x)g\big(\frac{x}{\mu_{B}(x)}\big)- \mu_{C}(x)\big|=\\
& &\big|\mu_{B}(x)g\big(\frac{x}{\mu_{B}(x)}\big)-
\mu_{B}(x)\mu_{C}\big(\frac{x}{\mu_{B}(x)}\big)\big|\leq
\mu_{B}(x)\varepsilon\big(\frac{x}{\mu_{B}(x)}\big)=\varepsilon_{0}\mu_{C}(x)
\end{eqnarray*}
for all $x\neq 0$. In particular, $\psi(x)\geq
(1-\varepsilon_{0})\mu_{C}(x)>0$ if $x\neq 0$. Therefore,
$$A:=\{x\in X : \psi(x)\leq 1\}$$ is a $C^p$ smooth starlike body
with respect to $0$. Let us check that $A$ approximates $C$ as
required. We have
\begin{eqnarray*}
& &x\in A\iff\psi(x)\leq 1\implies\mu_{C}(x)\leq
1+\varepsilon_{0}\mu_{C}(x)
    \implies (1-\varepsilon_{0})\mu_{C}(x)\leq 1\implies\\
    & &x\in\frac{1}{1-\varepsilon_{0}}C\subset (1+\delta)C,
\end{eqnarray*}
so $A\subset (1+\delta) C$. On the other hand, if $x\in
(1-\delta)C$, that is, $\mu_{C}(x)\leq 1-\delta$, then we have
    $$
    \psi(x)\leq
    (1+\varepsilon_{0})\mu_{C}(x)\leq(1+\varepsilon_{0})(1-\delta)<1,
    $$
hence $x\in A$.
\end{proof}

\begin{lem}\label{separation of K}
Let $X$ be a Banach space with $C^p$ smooth partitions of unity,
$K$ a weakly compact subset of $X$, and $D$ a bounded starlike
body with respect to $0$, such that $K\subset_{d} D$. Then there
exist $D_{1}$ and $D_{2}$, $C^p$ smooth starlike bodies with
respect to $0$, such that
    $$
    K\subset D_{1}\subset_{\mu} D_{2}\subset D.
    $$
\end{lem}
\begin{proof}
Since $K\subset_{d} D$ we can take $0<\theta<1/2$ so that
$K\subset (1-2\theta)D$. Choose $\delta\in (0,1)$ with
$(1-2\theta)/(1-\theta)<1-\delta$ and $(1+\delta)(1-\theta)<1$.
Applying the preceding lemma to $C:=(1-\theta)D$, we get a $C^p$
smooth starlike body with respect to $0$, $D_{1}$, such that
    $(1-\delta)C\subset D_{1}\subset (1+\delta)C$; in particular,
taking into account that $1-2\theta<(1-\theta)(1-\delta)$, we
deduce $K\subset (1-2\theta)D\subset
(1-\theta)(1-\delta)D=(1-\delta)C\subset D_{1}$. Now pick
$\varepsilon>0$ such that
$(1+\varepsilon)(1+\delta)(1-\theta)<1$, and set
$D_{2}:=(1+\varepsilon)D_{1}$. The body $D_{2}$ is $C^p$ smooth
and starlike with respect to $0$, and  $D_{1}\subset_{\mu}
D_{2}$. Finally, we also have
$D_{2}=(1+\varepsilon)D_{1}\subset(1+\varepsilon)(1+\delta) C
    \subset (1+\varepsilon)(1+\delta)(1-\theta) D\subset D$.
\end{proof}

\begin{lem}\label{interpolation of smooth bodies}
Let $X$ be a Banach space with $C^p$ smooth partitions of unity,
$C_{1}$ a bounded starlike body with respect to a point $c$, and
$C_{2}$ a mere subset of $X$ such that $C_{1}\subset_{d} C_{2}$.
Then there exist $D_{1}$ and $D_{2}$, $C^p$ smooth starlike bodies
with respect to $c$, which satisfy $C_{1}\subset_{\mu}
D_{1}\subset_{\mu} D_{2}\subset_{d} C_{2}$.
\end{lem}
\begin{proof}
We may assume that $c=0$ and $C_{2}\subset B_{X}$. Let us pick
$\varepsilon>0$ such that $\textrm{dist}(C_{1}, X\setminus
C_{2})\geq \varepsilon$. According to Lemma \ref{smooth
approximation of starlike bodies}, there exists a $C^p$ smooth
starlike body with respect to $0$, $A$, satisfying
    $$
    (1-\theta)(1+\varepsilon/2)C_{1}\subset A\subset (1+\theta)(1+\varepsilon/2)C_{1},
    $$
where $\theta$ is any positive number such that
$\varepsilon/2-\theta(1+\varepsilon/2)\geq\varepsilon/4$. Define
$D_{1}:=A$. Since $(1-\theta)(1+\varepsilon/2)\geq
1+\varepsilon/4$ we have that
$C_{1}\subset(1+\varepsilon/4)C_{1}\subset
(1-\theta)(1+\varepsilon/2)C_{1}\subset A=D_{1}$, and in
particular $C_{1}\subset_{\mu} D_{1}$. On the other hand,
    $$
    A=D_{1}\subset (1+\theta)(1+\varepsilon/2)C_{1}\subset
    C_{1}+(\theta+\varepsilon/2+\theta\varepsilon/2)C_{1}\subset
    C_{1}+(\theta+\varepsilon/2+\theta\varepsilon/2)B_{X},
    $$
which implies that
    $$
    \textrm{dist}(D_{1}, X\setminus C_{2})\geq
    \textrm{dist}\big(C_{1}+(\theta+\varepsilon/2+\theta\varepsilon/2)B_{X},
    X\setminus C_{2}\big)\geq\varepsilon-(\theta+\varepsilon/2+\theta\varepsilon/2)
    \geq\varepsilon/4.
    $$
Define now $D_{2}:=(1+\varepsilon/8)D_{1}$. It is obvious that
$D_{2}$ is a $C^p$ smooth starlike body with respect to $0$
satisfying $D_{1}\subset_{\mu} D_{2}$. Finally, we have that
$D_{2}=(1+\varepsilon/8)D_{1}\subset D_{1}+\varepsilon/8B_{X}$,
and therefore
    $$
    \textrm{dist}\big(D_{2}, X\setminus C_{2}\big)\geq
    \textrm{dist}\big(D_{1}+\frac{\varepsilon}{8}B_{X}, X\setminus C_{2}\big)
    \geq \textrm{dist}\big(D_{1}, X\setminus C_{2}\big)-\frac{\varepsilon}{8}
    \geq\frac{\varepsilon}{4}-\frac{\varepsilon}{8}=\frac{\varepsilon}{8},
    $$
which means that $D_{2}\subset_{d} C_{2}$.
\end{proof}

The following lemma is one of the keys to the proof of Proposition
\ref{existence of enough starlike bodies}.
\begin{lem}\label{tower of convex bodies with empty intersection}
Let $X$ be a nonreflexive Banach space, $K$ a weakly compact set,
and $C$ a bounded convex body with $0\in\textrm{int}(C)$ and
$K\subset_{d} C$. Then there exist $\varepsilon>0$ and a sequence
$(C_{n})$ of convex bodies such that
\begin{enumerate}
\item $\bigcap_{n=1}^{\infty}C_{n}=\emptyset$,
\item $C_{n+1}\subset_{d}C_{n}\subset C$ for all
$n\in\mathbb{N}$, and
\item $[K,C_{1}]+3\varepsilon B_{X}\subset C$.
\end{enumerate}
\end{lem}
\begin{proof}
Since $K\subset_{d} C$, there exists $\delta_{0}>0$ such that
$K\subset (1-2\delta_{0})C$ and, by Lemma \ref{smaller multiples
are at positive distance}, $\textrm{dist}\big((1-\delta_{0})C,
X\setminus C\big)\geq\delta_{1}$ for some $\delta_{1}>0$.

Since $X$ is nonreflexive, according to James' theorem, there
exists a continuous linear functional $T\in X^{*}$ such that $T$
does not attain its sup on the body $(1-2\delta_{0})C$,
    $$
    \alpha:=\sup\{ T(x) :x\in (1-2\delta_{0})C\}.
    $$
Define then $$H_{n}:=\{x\in (1-2\delta_{0})C : T(x)\geq
\alpha-1/n\}$$ for each $n\in\mathbb{N}$. We have that
    $$
    \bigcap_{n=1}^{\infty}H_{n}=\emptyset, \, \, H_{n+1}\subset H_{n}
    \textrm{ for all $n$, and }
    H_{1}\subset(1-2\delta_{0})C\subset_{d}(1-\delta_{0})C.
    $$
Take $\varepsilon>0$ such that $H_{1}+\varepsilon B_{X}\subset
(1-\delta_{0})C$ and $3\varepsilon<\delta_{1}$. Now, for each
$n\in\mathbb{N}$ let us define
    $$
    C_{n}=N(H_{n},\frac{\varepsilon}{2^{n}})=
    \{x\in X : \textrm{dist}(x, H_{n})\leq\frac{\varepsilon}{2^{n}}\}.
    $$
It is easy to see that $(C_{n})$ is a sequence of bounded convex
bodies such that
    $$
    \bigcap_{n=1}^{\infty}C_{n}=\emptyset, \textrm{ and }
    C_{n+1}\subset_{d} C_{n} \textrm{ for all } n\in\mathbb{N}.
    $$
Moreover, from the facts that
$C_{1}=N(H_{1},\varepsilon/2)\subset (1-\delta_{0})C$, $K\subset
(1-\delta_{0})C$, and
$3\varepsilon<\delta_{1}\leq\textrm{dist}((1-\delta_{0})C,
X\setminus C)$, we easily get that $[C_{1}, K]+3\varepsilon
B_{X}\subset C$.
\end{proof}

\medskip

\begin{center}
{\bf Proof of Proposition \ref{existence of enough starlike
bodies}.}
\end{center}

\noindent{\bf Case I.} Assume that $X$ is nonreflexive.

Let $E$ be a bounded convex body with $0\in\textrm{int}(E)$. By
Lemma \ref{smaller multiples are at positive distance}, we have
that $(1/2)E\subset_{d} E$. According to Lemma \ref{smooth
approximation of starlike bodies}, there exists a $C^p$ smooth
starlike body with respect to $0$ such that $(1/2)E\subset
B\subset E$. This body $B$ is the one we need.

Now take a weakly compact set $K\subset X$ and a starlike body
$A\supset K$. Since $A$ is bounded and $0\in\textrm{int}(E)$,
there exists some $r>0$ so that $A\subset_{d}(r/2)E$. According
to Lemma \ref{tower of convex bodies with empty intersection},
there exists $\varepsilon>0$ and a sequence $(C_{n})$ of convex
bodies such that $$ \bigcap_{n=1}^{\infty}C_{n}=\emptyset, \, \,
[K, C_{1}]+ 3\varepsilon B_{X}\subset (r/2)E, \, \textrm{ and }
C_{n+1}\subset_{d} C_{n}\subset rE \textrm{ for all }
n\in\mathbb{N}. $$ Let us choose a sequence $(c_{n})$ of points
of $X$ such that $c_{n}\in\textrm{int}(C_{n})$ for every
$n\in\mathbb{N}$. Set $\Delta=\textrm{diam}(\frac{r}{2}E)>0$. For
each $n\in\mathbb{N}$, define
    $$
    V_{n}=[C_{n}, K].
    $$
By Lemma \ref{that union of cones is starlike}, $V_{n}$ is a
Lipschitz starlike body with respect to every point in the
interior of $C_{n}$. Let $\mu_{n}=\mu_{V_{n}}$ be the Minkowski
functional of $V_{n}$ with respect to the point
$c_{n+1}\in\textrm{int}(C_{n+1})\subset\textrm{int}(C_{n})$. Note
that $\mu_{n}$ is a Lipschitz function.

Next we are going to inductively construct a sequence of positive
numbers $(\delta_{n})$ such that, if we define
    $$
    P_{n}:=\{x\in X : \mu_{n}(x)\leq 1+\delta_{n}\}
    $$
for each $n\in\mathbb{N}$, then $(P_{n})$ is a sequence of
bounded starlike bodies such that
\begin{itemize}
\item[{(i)}] $P_{n+1}\subset_{d} P_{n}\subset P_{1}\subset (r/2)E$
for all $n\in\mathbb{N}$,
\item[{(ii)}] $\bigcap_{n=1}^{\infty}P_{n}=K$,
\item[{(iii)}] $P_{n}$ is starlike with respect to $c_{n+1}$
for all $n\in\mathbb{N}$,
\item[{(iv)}] $C_{n+1}\subset_{d} P_{n}\cap C_{n}$ for all
$n\in\mathbb{N}$.
\end{itemize}

\noindent $\bullet${\em 1st step.} Choose $\delta_{1}>0$ with
$\delta_{1}<\min\{\varepsilon/\Delta, 1\}$, and set $P_{1}=\{x\in
X : \mu_{1}(x)\leq 1+\delta_{1}\}$. By Lemma \ref{Lipschitz
implies that mu implies d}, there is $\delta_{1}'>0$ such that
$P_{1}\supset V_{1}+\delta_{1}'B_{X}$.

\noindent $\bullet${\em 2nd step.} Now choose $\delta_{2}>0$ such
that $\delta_{2}<\min\{\delta_{1}'/2\Delta, 1/2\}$. Then
$P_{2}=\{x\in X : \mu_{2}(x)\leq 1+\delta_{2}\}\subset
V_{2}+(\delta_{1}'/2)B_{X}$, and therefore $\textrm{dist}(P_{2},
X\setminus P_{1})>0$.

\noindent $\bullet${\em (n+1)-th step.} Assume $\delta_{j}$ and
$P_{j}$ are already defined for $j=1, 2, ..., n$ in such a way
that $P_{j+1}\subset_{d} P_{j}$ for $j\leq n-1$. By Lemma
\ref{Lipschitz implies that mu implies d}, there is
$\delta_{n}'>0$ such that $P_{n}\supset V_{n}+\delta_{n}'B_{X}$.
Pick $\delta_{n+1}>0$ so that
$\delta_{n+1}<\min\{\delta_{n}'/2\Delta, 1/2^{n}\}$, and set
$P_{n+1}=\{x\in X : \mu_{n+1}(x)\leq 1+\delta_{n+1}\}$. Then we
have that $P_{n+1}\subset V_{n}+(\delta_{n}'/2)B_{X}$, hence
$\textrm{dist}(P_{n+1}, X\setminus P_{n})>0$.

\medskip

By induction the sequence $(P_{n})$ is well-defined and satisfies
properties $(i)$ and $(iii)$ above. To see that $P_{1}\subset_{d}
(r/2)E$, just note that $P_{1}\subset V_{1}+\delta_{1}\Delta
B_{X}=[C_{1},K]+\delta_{1}\Delta B_{X}\subset [C_{1},
K]+3\varepsilon B_{X}\subset (r/2)E$. On the other hand, since
$P_{n}\cap C_{n}=C_{n}$, it is clear that $C_{n+1}\subset_{d}
P_{n}\cap C_{n}$, that is, the sequence $(P_{n})$ satisfies
property $(iv)$.

Finally, let us check that condition $(ii)$ is met as well. It is
immediate that $K\subset\bigcap_{n=1}^{\infty}P_{n}$. Let us take
$q\in\bigcap_{n=1}^{\infty}P_{n}$ and show that $q\in K$. For
each $n\in\mathbb{N}$ we have $q\in P_{n}\subset
V_{n}+\delta_{n}\Delta B_{X}=[C_{n},K]+\delta_{n}\Delta B_{X}$,
so there are $x_{n}\in C_{n}, y_{n}\in K, t_{n}\in [0,1]$ with
$\|q-(1-t_{n})x_{n}-t_{n}y_{n}\|\leq\delta_{n}\Delta$, and in
particular $\lim_{n\to\infty}[(1-t_{n})x_{n}-t_{n}y_{n}]=q$. Since
$K$ is weakly compact and $[0,1]$ is compact, we may assume
(passing to a subsequence if necessary) that $y_{n}$ converges to
some $y_{0}\in K$ weakly, and $t_{n}\to t_{0}\in [0,1]$. Then
$(1-t_{n})x_{n}$ converges to $q-t_{0}y_{0}$ weakly. If
$t_{0}\neq 1$ then we have that $x_{n}$ converges weakly to
$x_{0}:=(1-t_{0})^{-1}(q-t_{0}y_{0})$; but, since each
$C_{n}\supset (x_{j})_{j\geq n}$ is closed and convex, hence
weakly closed, we have $x_{0}\in C_{n}$ for each $n$, and then
$x_{0}\in\bigcap_{n=1}^{\infty}C_{n}=\emptyset$, a contradiction.
Therefore, $t_{0}=1$, and $q=y_{0}\in K$.

\smallskip

Now we are going to define the bodies $A_{n}, B_{n}, D_{n},
E_{n}$, and $Q_{n}$. Fix $n\in\mathbb{N}$. Since $C_{n+2}$ and
$C_{n+1}$ are bounded starlike bodies with respect to $c_{n+2}$,
and $C_{n+2}\subset_{d} C_{n+1}$, we can apply Lemma
\ref{interpolation of smooth bodies} to obtain two $C^p$ smooth
starlike bodies $D_{n}, E_{n}$ with respect to $c_{n+2}$ such that
    $$
    C_{n+2}\subset D_{n}\subset_{\mu} E_{n}\subset_{d} C_{n+1}.
    $$
Another application of Lemma \ref{interpolation of smooth bodies}
gives us a $C^p$ smooth starlike body $A_{n}$ with respect to
$c_{n+2}$ such that
    $$
    E_{n}\subset_{\mu} A_{n}\subset C_{n+1}=C_{n+1}\cap P_{n+1}.
    $$
Besides, $P_{n+1}\subset_{d} P_{n}$, and $P_{n+1}$ is starlike
with respect to $c_{n+1}$. Then, applying Lemma
\ref{interpolation of smooth bodies} for the last time (now $P_n$
acts as a mere set, it is not necessary that $P_n$ be starlike
with respect to $c_{n+2}$, only $P_{n+1}$ has to meet this
condition), we get $B_{n}$ and $Q_{n}$, two $C^p$ smooth starlike
bodies with respect to $c_{n+2}$, satisfying
    $$
    P_{n+1}\subset_{\mu} B_{n}\subset_{\mu} Q_{n}\subset P_{n}.
    $$
Moreover, we also have $E_{n}\subset C_{n+1}\subset
P_{n+1}\subset_{\mu} B_{n}$. Summing up, we get that
    $$
    C_{n+2}\subset D_{n}\subset_{\mu} E_{n}\subset_{\mu} A_{n}
    \subset C_{n+1}\subset P_{n+1}\subset_{\mu} B_{n}\subset_{\mu}
    Q_{n}\subset P_{n},
    $$
and now it is clear that the sequences of bodies we have just
constructed satisfy conditions $(1)-(4)$ of Proposition
\ref{existence of enough starlike bodies}. Finally, $B$ is the
required universal body and satisfies condition $(5)$. Indeed,
notice that $K\subset A\subset (r/2)\textrm{int}(E)\subset rB$,
$P_{1}\subset (r/2)E\subset rB$.

\medskip

\noindent{\bf Case II.} Assume now that $X$ is reflexive.

In this case it is known that there exists a continuous linear
injection $T:X\longrightarrow c_{0}(\Gamma)$ for some (infinite)
set $\Gamma$ (see \cite{DGZ}, p.246, for instance). It is also
well known that for an infinite set $\Gamma$, the space
$c_{0}(\Gamma)$ is $c_{0}$-saturated, that is, every
infinite-dimensional closed subspace of $c_{0}(\Gamma)$ has a
closed subspace which is isomorphic to $c_{0}$. This clearly
implies that $c_{0}(\Gamma)$ contains no closed
infinite-dimensional reflexive subspaces. Therefore
$Y:=\overline{T(X)}\subset c_{0}(\Gamma)$ is nonreflexive, and
$T(X)$ is not a closed subspace of $Y\subset c_{0}(\Gamma)$. On
the other hand, the space $c_{0}(\Gamma)$ has a $C^{\infty}$
smooth equivalent norm (see \cite{DGZ}, chapter V, theorem 1.5),
whose restriction to $Y$ defines a $C^\infty$ smooth equivalent
norm $|\cdot|$. Finally, it is well known \cite{DGZ} that the
space $c_{0}(\Gamma)$ has $C^\infty$ smooth partitions of unity,
hence so does $Y$.

Summing up, we have a continuous linear injection $T:X\to Y$,
where $(Y, |\cdot|)$ is a nonreflexive Banach space with a
$C^\infty$ smooth norm and $C^\infty$ smooth partitions of unity,
and $T(X)$ is dense in $Y$.

Set $B'=\{y\in Y : |y|\leq 1\}$, which is a $C^\infty$ smooth
bounded convex body with $0\in\textrm{int}(B')$. Define
$B=T^{-1}(B')$. It is clear that $B$ is a radially bounded
$C^\infty$ smooth convex body.

Let $K$ be a weakly compact subset of $X$ and $A$ a bounded
starlike body containing $K$. Since $T$ is continuous, $T(A)$ is
bounded in $Y$. Choose $r>0$ so that $T(A)\subset_{d}(r/2)B'$.
Now we may copy the above proof (nonreflexive case) with $B'=E$
and obtain sequences of $C^\infty$ smooth starlike bodies,
$(P_{n}'), (C_{n}'), (A_{n}'), (B_{n}'), (Q_{n}'), (D_{n}'),
(E_{n}')$, and a sequence of points $(c_{n}')$ of $Y$ satisfying
the conditions $(1)-(4)$ of the statement of Proposition
\ref{existence of enough starlike bodies} and $P_{1}'\subset
(r/2)B'$. Ensure further that $c_{n}'\in
T(X)\cap\textrm{int}(C_{n}')$ for each $n\in\mathbb{N}$ (this is
possible because $T(X)$ is dense in $Y$, hence
$T(X)\cap\textrm{int}(C_{n}')\neq\emptyset$ for all $n$).

Then, for each $n\in\mathbb{N}$, define $c_{n}=T^{-1}(c_{n}')\in
X$, and
\begin{eqnarray*}
& & C_{n}=T^{-1}(C_{n}'), \, B_{n}=T^{-1}(B_{n}'), \,
P_{n}=T^{-1}(P_{n}'), \, A_{n}=T^{-1}(A_{n}')\\&
&Q_{n}=T^{-1}(Q_{n}'), \,  D_{n}=T^{-1}(D_{n}'), \,
E_{n}=T^{-1}(E_{n}')\subset X.
\end{eqnarray*}
By Lemma \ref{pullback}, these are radially bounded $C^\infty$
smooth starlike bodies with respect to $c_{n+2}$. On the other
hand, Lemma \ref{pullback preserves strict inclusions} guarantees
that $$
    C_{n+2}\subset D_{n}\subset_{\mu} E_{n}\subset_{\mu} A_{n}
    \subset C_{n+1}\subset P_{n+1}\subset B_{n}\subset_{\mu}
    Q_{n}\subset P_{n}.
     $$
Finally, it is immediately checked that
$\bigcap_{n=1}^{\infty}C_{n}=\emptyset$,
$\bigcap_{n=1}^{\infty}P_{n}=K$, $P_{1}\subset (r/2)B$, and
$A\subset (r/2)B$.

\medskip

\noindent Now we are in a position to finish the proof of the main
result.

\medskip

\begin{center}
{\bf Proof of Theorem \ref{main thm}.}
\end{center}

\noindent We may assume that $A$ is starlike with respect to the
origin. According to Proposition \ref{existence of enough
starlike bodies}, there exists a radially bounded $C^p$ smooth
starlike body $B$ with respect to $0$ so that $K\subset_{d}
A\subset_{\mu} rB$ for some $r>0$ large enough, and there are
sequences $(P_{n})$, $(C_{n})$, $(A_{n})$, $(B_{n})$, $(Q_{n})$,
$(D_{n})$, $(E_{n})$ of subsets of $X$ and a sequence $(c_{n})$
of points of $X$ which satisfy the conditions of Proposition
\ref{construction of the deleting diffeomorphism}. Then we can
apply this Proposition to find a $C^p$ diffeomorphism $\Psi:X\to
X\setminus K$ such that $\Psi$ is the identity on $X\setminus
P_{1}\supset X\setminus rB$.

On the other hand, since $K\subset_{d} A$, Lemma \ref{separation
of K} allows us to find two $C^p$ smooth starlike bodies $U_{1},
U_{2}$ with respect to $0$ such that
    $$
    K\subset U_{1}\subset_{\mu} U_{2}\subset A.
    $$
Now, by the Four Bodies Lemma \ref{four bodies}, there is a $C^p$
diffeomorphism $g:X\to X$ such that $g(U_{2})=rB$ and $g$ is the
identity on $U_{1}\supset K$; notice in particular that $g(K)=K$.

Define then $h=g^{-1}\circ\Psi\circ g$. It is clear that $h$ is a
$C^p$ diffeomorphism from $X$ onto $X\setminus K$. Moreover, if
$x\in X\setminus A$ then $x\notin U_{2}$, so $g(x)\notin rB$,
which implies that $\Psi(g(x))=g(x)$, hence $h(x)=x$; that is, $h$
is the identity off $A$. 

\medskip

\begin{center}
{\bf Acknowledgements}
\end{center}
\noindent Part of this research was carried out during a stay of
the first-named author in the Mathematics Department of University
College London. This author wishes to thank the members of that
Department, and very especially David Preiss, for their kind
hospitality and advice. The authors also wish to thank Antonio
Su\'arez Granero for his probity.



\vspace{3mm} \noindent Departamento de An\'alisis Matem\'atico.
Facultad de Ciencias Matem\'aticas. Universidad Complutense.
28040 Madrid, SPAIN\\  \noindent {\em E-mail addresses:}
daniel\_azagra@mat.ucm.es, a\_montesinos@mat.ucm.es
\end{document}